\documentclass[12pt,amsthm,seceqn]{elsart}

\usepackage{amssymb,amsmath}
\usepackage[all]{xy}

\newcommand{\vp}{\varepsilon}
\newcommand{\cl}[1]{{\mathcal{#1}}}
\newcommand{\bb}[1]{{\mathbb{#1}}}

\newcommand{\AcG}{\mathcal{A}\times_\alpha G}
\newcommand{\coA}{(\AcG)\times_{\hat\alpha}G}
\newcommand{\AK}{\mathcal{A}\otimes\mathcal{K}(L^2(G))}
\newcommand{\adlg}{Ad_{\lambda_g}}
\newcommand{\adl}{Ad_{l_gr_g}}
\newcommand{\adlrg}{Ad_{\lambda_g\otimes\l_gr_g}}

\theoremstyle{plain}
\newtheorem{lemma}{Lemma}[section]
\newtheorem{proposition}[lemma]{Proposition}
\newtheorem{mythm}[lemma]{Theorem}

\theoremstyle{definition}

\theoremstyle{remark}
\newtheorem{myrem}[lemma]{Remark}

\begin{document}

\begin{frontmatter}

  \title{Crossed Products and Entropy of Automorphisms} \author{Ciprian
    Pop\corauthref{current}}
  \address{Institute of Mathematics of the Romanian Academy\\
    C.P. 1--764\\ Bucharest, Romania} \corauth[current]{Current Adress:
    Department of Mathematics Texas A\&M University College Station, TX
    77843} \ead{cpop@math.tamu.edu} \author{Roger R.
    Smith\thanksref{nsf}} \ead{rsmith@math.tamu.edu} \address{Department of
    Mathematics\\ Texas A\&M University\\ College Station, TX 77843--3368}
  \thanks[nsf]{Partially supported by a grant from the National Science
    Foundation.}

\begin{abstract}
  Let $\cl A$ be an exact $C^*$--algebra, let $G$ be a locally compact
  group, and let $(\cl A, G, \alpha)$ be a $C^*$--dynamical system. Each
  automorphism $\alpha_g$ induces a spatial automorphism $Ad_{\lambda_g}$
  on the reduced crossed product $\AcG$. In this paper we examine the
  question, first raised by E.~St{\o}rmer, of when the topological
  entropies of $\alpha_g$ and $Ad_{\lambda_g}$ coincide. This had been
  answered by N.~Brown for the particular case of discrete abelian groups.
  Using different methods, we extend his result to preservation of entropy
  for $\alpha_g$ when the subgroup of ${\text{Aut}}(G)$ generated by the
  corresponding inner automorphism $Ad_g$ has compact closure. This
  property is satisfied by all elements of a wide class of groups called
  locally $[FIA]^-$. This class includes all abelian groups, both discrete
  and continuous, as well as all compact groups.
\end{abstract}

\end{frontmatter}

\newpage
\section{Introduction}\label{sec1}

\indent

In \cite{V}, Voiculescu introduced the topological entropy $ht(\alpha)$ of
an automorphism $\alpha$ of a nuclear $C^*$--algebra, generalizing the
classical notion for abelian $C^*$--algebras. The definition, which we
recall in Section \ref{sec2}, was based on the theorem of Choi and Effros,
\cite{CE}, which characterized nuclear $C^*$--algebras as those which
admit approximate point norm completely positive factorizations of the
identity map through matrix algebras. Brown observed in \cite{Br} that the
definition could be extended to exact $C^*$--algebras $\cl A$, by allowing
the completely positive approximations to have range in any containing
$B(H)$ rather than $\cl A$ itself. This was based on two results.
Wassermann, \cite{W}, had shown that the existence of such factorizations
implies exactness and subsequently Kirchberg, \cite{Ki}, proved that this
property characterizes exact $C^*$--algebras. Any automorphism $\alpha$ of
$\cl A$ induces an action of $\bb Z$ on $\cl A$, and then $\alpha$ is
implemented by an inner automorphism $Ad_u$ on the crossed product
${\cl A}\times_\alpha \bb Z$. In \cite{Sto}, St{\o}rmer had posed the
question of whether passing to the crossed product preserves the entropy,
the point being to replace an arbitrary automorphism by an inner one which
would be more amenable to analysis. Brown, \cite{Br}, answered this
positively for exact $C^*$--algebras. Since exactness is preserved by
crossed products of amenable groups, \cite{Ki2}, St{\o}rmer's question is
relevant and interesting in this wider context. In this paper we show that
an affirmative answer can be given for a large class of locally compact
amenable groups.

Brown's approach, stated for $\bb Z$ but valid for any discrete abelian
group, was based on work in \cite{SS}, where completely {\it{positive}}
factorizations of $\cl A$ through $\AcG$ were constructed for discrete
groups. We do not know whether such factorizations exist beyond the discrete 
case,
but we were able to find completely {\it{contractive}} factorizations of $\cl 
A$
through $\AcG$ for general locally compact groups, \cite{NS}. In order to make 
use of this result,
we have been led to introduce a new entropy $ht'(\alpha)$ for an automorphism
$\alpha$ involving
complete contractions rather than completely positive maps. The first goal of 
the 
paper is to show that $ht'(\alpha)=ht(\alpha)$ (see Theorem \ref{thm3.7}), and 
then the 
point of defining $ht'(\alpha)$ is to have an equivalent formulation of
$ht(\alpha)$ which is more suited to our context. The
proof of equality of these two entropies relies, in part, on techniques
from \cite{Sm}.
 
The second section reviews some background material and the third
discusses our new definition of entropy. The results of Section \ref{sec5}
ultimately depend on non--abelian duality as developed in \cite{IT,T}, and
we state a modified version of the Imai--Takai duality theorem in Theorem
\ref{thm4.2}, essentially due to Quigg, \cite{Q}, in order to obtain extra 
information not available from
\cite{IT}. Section \ref{sec4} also shows how entropy changes as
automorphisms are lifted from $\cl A$ to the crossed product $\AcG$ and
then to the double crossed product $\coA$, using the theory of coactions.

When ${\mathcal A}={\mathbb C}$, the action is trivial and we have
$ht(\alpha_g)=0$. The crossed product becomes $C^*_r(G)$ and $\lambda_g$
is the lifting of the corresponding inner automorphism of $G$ to this
$C^*$--algebra. If entropy is to be preserved, then these automorphisms of
$C^*_r(G)$ must have zero entropy. This is clearly true for abelian groups
where the set of inner automorphisms reduces to a single point. A natural
extension of this case is to impose a compactness requirement on the closure
in ${\mathrm{Aut}}(G)$
of the subgroup  generated by the inner automorphism $Ad_g$ for a single group 
element
$g\in G$. This leads us to consider some classes of groups which are
standard in the literature, \cite{Pal}. Our
results are proved in the context of [{\it {SIN}}] groups and $[FIA]^-$ groups
(see \cite{Pal} for a detailed survey of classes of locally compact
groups). The first class is defined by the property of having small
invariant neighborhoods of the identity element $e$ under the inner 
automorphisms, while the
second class requires the closure of the set of inner automorphisms to be
compact in ${\text {Aut}}(G)$ for its natural topology (see \cite{Pet} for a
nice characterization of this). Groups in the second class are amenable and 
both classes contain all
abelian and all compact groups. In fact our results are valid for a much
larger class of groups which are amenable and satisfy a condition 
which we call locally $[FIA]^-$. The definition only
requires that the closure of the group generated by a single but arbitrary
inner automorphism should be compact in ${\text {Aut}}(G)$. These matters
are discussed in Section \ref{sec5} which contains our main results.

Throughout the paper we restrict $C^*$--algebras to be exact since, as
pointed out by Brown, \cite{Br}, topological entropy only makes sense for
this class. However, an exactness hypothesis is only necessary in Section 
\ref{sec2},
Theorems \ref{thm3.5} and \ref{thm3.7}, and Section \ref{sec5}. All other 
results are valid in general.

We thank John Quigg for helpful correspondence concerning Theorem \ref{thm4.2},
and also the referee for some enlightening comments.

\newpage

\section{Preliminaries}\label{sec2}

\indent

In this section we give a brief review of crossed products and entropy of 
automorphisms. We will recall the completely contractive factorization 
property 
(CCFP) from \cite{NS}, and define a covariant version (CCCFP) which will be 
useful subsequently. We will also introduce a new version of topological 
entropy, using completely contractive maps rather than completely positive  
ones. We will show that this version coincides with the standard definition, 
and 
the point is to obtain a more flexible situation in which to work.

Let $G$ be a locally compact group, let $\cl A$ be a $C^*$-algebra, and let 
$\alpha \colon \ G\to \text{Aut}(\cl A)$ be a homomorphism such that the map 
$t\mapsto \alpha_t(a)$ is norm continuous on $G$ for each $a\in\cl A$. Then 
the 
triple $(\cl A,G,\alpha)$ is called a $C^*$-dynamical system. The reduced 
crossed $\cl A \times_{\alpha,r} G$ is constructed by taking a faithful 
representation 
$\pi\colon \ \cl A\to B(H)$ and then defining a new representation $\tilde\pi 
\colon \ \cl A\to B(L^2(G,H))$ by
\begin{equation}\label{eq2.1}
(\tilde\pi(a)\xi)(t) = \alpha_{t^{-1}}(a)(\xi(t)),\qquad a\in\cl A, \ \xi\in 
L^2(G,H).
\end{equation}
There is an associated unitary representation of the group given by
\begin{equation}\label{eq2.2}
(\lambda_s\xi)(t) = \xi(s^{-1}t),\qquad s\in G,\ \xi \in L^2(G,H),
\end{equation}
and in this representation the action $\alpha$ is spatially implemented, as 
can 
be seen from the easily established equation
\begin{equation}\label{eq2.3}
\lambda_s\tilde\pi(a) \lambda_{s^{-1}} = \tilde\pi(\alpha_s(a)),\qquad a\in\cl 
A, \ s\in G.
\end{equation}
The norm closed span of operators on $L^2(G,H)$ of the form
\begin{equation}\label{eq2.4}
\int_G f(s) \tilde\pi(a)\lambda_s\ ds, \qquad a\in\cl A, \ f\in C_c(G),
\end{equation}
for a fixed left Haar measure $ds$, is a $C^*$-algebra called the reduced 
crossed product. There is also a full crossed product, denoted $\cl 
A\times_\alpha G$, which will not appear here. The groups considered will 
eventually be amenable, where the two notions coincide. For this reason we 
will 
abuse the standard notation and write $\cl A\times_\alpha G$ for the reduced 
crossed  product throughout. The construction is independent of the 
representation $\pi$, and so we may carry it out twice, if necessary, so that 
we 
may always assume that the action is spatially implemented by a unitary 
representation of $G$.

Let $\cl A$ and $\cl B$ be two $C^*$-algebras. If there exist nets of 
complete contractions
\begin{equation}\label{eq2.5}
\cl A \stackrel{S_\lambda}{\longrightarrow} \cl B 
\stackrel{T_\lambda}{\longrightarrow} \cl A
\end{equation}
so that
\begin{equation}\label{eq2.6}
\lim_\lambda\|T_\lambda(S_\lambda(a)) - a\| = 0,\qquad a\in\cl A,
\end{equation}
then we say that the pair $(\cl A,\cl B)$ has the completely contractive 
factorization property (CCFP), \cite{NS}. We will also be considering a pair 
of 
automorphisms $\alpha$ and $\beta$ of $\cl A$ and $\cl B$  respectively. If 
the 
above nets can be chosen to satisfy the additional requirement that
\begin{equation}\label{eq2.7}
S_\lambda\alpha = \beta S_\lambda,\quad T_\lambda\beta = \alpha T_\lambda,
\end{equation}
then the pair $(\cl A,\alpha), (\cl B,\beta)$ is said to have the covariant 
CCFP 
 (CCCFP).

In \cite{V}, Voiculescu introduced the notion of topological entropy for an 
automorphism of a nuclear $C^*$-algebra, and this was extended by Brown to the 
case of an exact $C^*$-algebra, \cite{Br}. We assume that $\cl A$ is 
concretely 
represented on a Hilbert space $H$, and we fix an automorphism $\alpha$ of 
$\cl 
A$. Given a finite subset $\omega\subseteq\cl A$ and $\delta > 0$, 
$rcc(\omega,\delta)$ denotes the smallest integer $m$ for which complete 
 contractions
\begin{equation}\label{eq2.8}
\cl A \stackrel{\phi}{\longrightarrow} \bb M_m 
\stackrel{\psi}{\longrightarrow} 
B(H)
\end{equation}
can be found to satisfy
\begin{equation}\label{eq2.9}
\|\psi(\phi(a)) - a\| < \delta,\qquad a\in\omega.
\end{equation}
We then define successively
\begin{align}
\label{eq2.10}
ht'(\alpha,\omega,\delta) &= \limsup_{n\to\infty} \frac1n \log \left(rcc\left( 
\bigcup^{n-1}_{i=0} \alpha^i(\omega),\delta\right)\right),\\
\label{eq2.11}
ht'(\alpha,\omega) &= \sup_{\delta>0} ht'(\alpha,\omega,\delta), \\
\label{eq2.12}
ht'(\alpha) &= \sup_\omega ht'(\alpha,\omega).
\end{align}
This last quantity is our version of the topological entropy $ht(\alpha)$ of 
$\alpha$,
originally defined as above using completely positive contractions throughout. 
The term
$rcc(\cdot)$ (derived from {\it{completely contractive rank}}) replaces the 
{\it{completely positive
rank}} $rcp(\cdot)$. Since positivity of maps is a more stringent requirement, 
the entropy 
arising from the completely contractive framework is immediately seen to 
satisfy  
\begin{equation}\label{eq2.13}
ht'(\alpha) \le ht(\alpha).
\end{equation}
We will  prove equality in the next section, so $ht'(\alpha)$ is a temporary 
notation soon to be replaced by $ht(\alpha)$. When confusion might arise, we 
will use
the notation $ht_{\cl A}(\alpha)$ to indicate the algebra on which $\alpha$ 
acts.
If $\alpha \in {\text {Aut}}(\cl A)$ and $\cl B$ is an $\alpha$--invariant 
$C^*$--subalgebra,
it may be the case that $ht_{\cl A}(\alpha)\ne ht_{\cl B}(\alpha)$.

The following lemmas will be useful subsequently in comparing the entropies of 
two automorphisms.

\begin{lemma}\label{lem2.1}
Let $\alpha$ and $\beta$ be automorphisms of exact $C^*$-algebras $\cl A$ and 
$\cl B$ respectively.
\begin{itemize}
\item[(i)] Suppose that, given an arbitrary finite subset $\omega\subseteq \cl 
A$ and 
$\delta>0$, there exist a finite subset $\omega'\subseteq \cl B$ and 
$\delta'>0$ 
so that
\begin{equation}\label{eq2.14}
ht'_{\cl A}(\alpha,\omega,\delta) \le ht'_{\cl B}(\beta,\omega',\delta').
\end{equation}
Then $ht'_{\cl A}(\alpha) \le ht'_{\cl B}(\beta)$.
\item[(ii)] Suppose that there exists a function $g\colon \ (0,\infty)\to 
(0,\infty)$ satisfying $\lim\limits_{\delta\to 0} g(\delta) = 0$, and suppose 
that, given an arbitrary finite subset $\omega\subseteq \cl A$ and $\delta>0$, 
there exists 
a finite subset $\omega'\subseteq \cl B$ such that
\begin{equation}\label{eq2.15}
ht'_{\cl A}(\alpha,\omega,g(\delta)) \le ht'_{\cl B}(\beta,\omega',\delta).
\end{equation}
Then $ht'_{\cl A}(\alpha) \le ht'_{\cl B}(\beta)$.
\end{itemize}
\end{lemma}

\begin{proof}
These are simple consequences of the definitions, after observing that 
$ht'(\alpha,\omega,\delta)$ is a monotonically decreasing function of $\delta$ 
for fixed $\alpha$ and $\omega$.
\end{proof}

\begin{lemma}\label{lem2.2}
Let $\alpha$ and $\beta$ be respectively automorphisms of the exact 
$C^*$-algebras $\cl A\subseteq B(H)$ and $\cl B\subseteq B(K)$, and suppose 
that 
the pair $(\cl A,\alpha)$, $(\cl B,\beta)$ has the CCCFP. Then $ht'_{\cl 
A}(\alpha) \le 
ht'_{\cl B}(\beta)$.
\end{lemma}

\begin{proof}
Let
\begin{equation}\label{eq2.16}
\cl A \stackrel{S_\lambda}{\longrightarrow} \cl B 
\stackrel{T_\lambda}{\longrightarrow} \cl A
\end{equation}
be nets of complete contractions satisfying \eqref{eq2.7}. Fix a finite subset 
$\omega\subseteq \cl A$, and $\delta>0$. Choose $\lambda$ so that
\begin{equation}\label{eq2.17}
\|T_\lambda(S_\lambda(a)) - a\| < \delta,\qquad a\in\omega,
\end{equation}
and let $\omega' = S_\lambda(\omega) \subseteq \cl B$. For each integer $n$, 
let 
$r_n$ be $rcc\left(\bigcup\limits^{n-1}_{i=0} \beta^i(\omega'), 
\delta\right)$. 
Then, by definition, there exist complete contractions
\begin{equation}\label{eq2.18}
\cl B \stackrel{\phi_n}{\longrightarrow} \bb M_{r_n} 
\stackrel{\psi_n}{\longrightarrow} B(K)
\end{equation}
such that
\begin{equation}\label{eq2.19}
\|\psi_n\phi_n(b) - b\| < \delta, \qquad b\in \beta^i(\omega'), \ 0 \le i \le 
n-1.
\end{equation}
By injectivity of $B(H)$, $T_\lambda$ has a completely contractive extension 
$\widetilde 
T_\lambda\colon \ B(K)\to B(H)$. Consider the diagram
\begin{equation}\label{eq2.20}
\cl A \stackrel{S_\lambda}{\longrightarrow} \cl B 
\stackrel{\phi_n}{\longrightarrow} \bb M_{r_n} 
\stackrel{\psi_n}{\longrightarrow} B(K) \stackrel{\widetilde 
T_\lambda}{\longrightarrow} B(H).
\end{equation}
If $a\in\omega$ then
\begin{equation}\label{eq2.21}
S_\lambda \alpha^i(a) = \beta^i(S_\lambda(a)),\quad T_\lambda\beta^i 
(S_\lambda(a)) = \alpha^i(T_\lambda(S_\lambda(a))),
\end{equation}
and a simple approximation argument shows that 
$rcc\left(\bigcup\limits^{n-1}_{i=0} \alpha^i(\omega), 2\delta\right) \le 
r_n$. 
It follows that
\begin{equation}\label{eq2.22}
ht'(\alpha,\omega,2\delta) \le ht'(\beta,\omega',\delta).
\end{equation}
The result follows from Lemma \ref{lem2.1} (ii) with $g(\delta) = 2\delta$.
\end{proof}
\newpage

\section{Completely contractive topological entropy}\label{sec3}

\indent

The objective in this section is to show that the two entropies $ht(\alpha)$ 
and 
$ht'(\alpha)$ coincide. This will be accomplished in two stages. The first is 
to 
show preservation of $ht'(\cdot)$ when lifting an automorphism from $\cl A$ to 
its unitization $\tilde{\cl A}$. The second is to show equality for 
automorphisms of unital $C^*$-algebras. Throughout $\cl A$ is assumed to be 
exact, and so $\tilde{\cl A}$ is also exact, \cite{Ki}. Various norm estimates 
and approximations will be needed and so we begin with some technical lemmas.
Many of the complications stem from allowing $\cl A$ to be non--unital. These 
are unavoidable
since crossed products inevitably take us out of the category of unital 
$C^*$--algebras. 
Below, $\cl A^+_1$ denotes the positive part of the closed unit ball.

\begin{lemma}\label{lem3.1}
Let $a,p\in \cl A^+_1$ satisfy $\|a-ap\|\le C$ for some $C>0$. If  $0\le b\le 
a$, 
then $\|b-bp\| \le\sqrt C$.
\end{lemma}

\begin{proof}
This follows from the inequality
\begin{equation}\label{eq3.1}
0 \le (1-p)b^2(1-p) \le (1-p) b(1-p) \le (1-p)a(1-p) \le C1,
\end{equation}
which gives $\|b(1-p)\|^2 \le C$.
\end{proof}

\begin{lemma}\label{lem3.2}
If $t\in B(H)$ satisfies
\begin{equation}\label{eq3.2}
\|1+e^{i\theta}t\| \le 1+C,\qquad\theta\in \bb R,
\end{equation}
for some constant $C>0$, then $\|t\| \le 2C$.
\end{lemma}

\begin{proof}
If $\phi$ is a state on $B(H)$, then
\begin{equation}\label{eq3.3}
|1+e^{i\theta}\phi(t)| \le 1+C,\qquad \theta\in\bb R.
\end{equation}
A suitable choice of $\theta$ shows that $|\phi(t)|\le C$, so the numerical 
radius is at most $C$, and the result follows.
\end{proof}

\begin{lemma}\label{lem3.3}
Let $\phi\colon \ \cl A\to B(H)$ be a complete contraction on a unital 
$C^*$-algebra, and let $t$ be the self-adjoint part of $\phi(1)$. Then there 
exists a completely positive contraction $\psi\colon \ \cl A\to B(H)$ such that
\begin{equation}\label{eq3.4}
\|\phi-\psi\|_{cb} \le \sqrt{2\|1-t\|}.
\end{equation}
\end{lemma}

\begin{proof}
 From the representation theorem for completely bounded maps, \cite{Pau}, 
there 
exist a representation $\pi\colon\ \cl A\to B(K)$ and contractions 
$V_1,V_2\colon \ H\to K$ such that
\begin{equation}\label{eq3.5}
\phi(x) = V^*_1\pi(x) V_2,\qquad x\in\cl A.
\end{equation}
Then
\begin{align}
0 \le (V_1-V_2)^* (V_1-V_2) &= V^*_1V_1 + V^*_2V_2 - V^*_1V_2 - 
V^*_2V_1\nonumber\\
\label{eq3.6}&\le 2-2t
\le 2\|1-t\|,
\end{align}
so $\|V_1-V_2\| \le \sqrt{2\|1-t\|}$. If we define $\psi$ to be 
$V^*_1\pi(\cdot)V_1$, then \eqref{eq3.4} is immediate.
\end{proof}

\begin{lemma}\label{lem3.4}
Let $\omega$ be a subset of $\cl A_{s.a.}$, let $\vp>0$, and let
\begin{equation}\label{eq3.7}
\cl A \stackrel{\phi}{\longrightarrow} \bb M_n 
\stackrel{\psi}{\longrightarrow} 
B(H)
\end{equation}
be a diagram of complete contractions satisfying
\begin{equation}\label{eq3.8}
\|\psi(\phi(x)) - x\| < \vp, \qquad x\in\omega.
\end{equation}
Then there exist self-adjoint complete contractions
\begin{equation}\label{eq3.9}
\cl A \stackrel{\tilde\phi}{\longrightarrow} \bb M_{2n} 
\stackrel{\widetilde\psi}{\longrightarrow} B(H)
\end{equation}
satisfying
\begin{equation}\label{eq3.10}
\|\widetilde\psi(\tilde\phi(x)) - x\| < \vp,\qquad x\in\omega.
\end{equation}
\end{lemma}

\begin{proof}
Following the methods of Paulsen, \cite{Pau}, we define $\tilde\phi\colon\ \cl 
A\to \bb M_{2n}$ by
\begin{equation}\label{eq3.11}
\tilde\phi(a) = \left(\begin{matrix} 0&\phi(a)\\ \phi(a^*)^*&0\end{matrix} 
\right),\quad a\in\cl A,
\end{equation}
which is a self-adjoint complete contraction. The map
\begin{equation}\label{eq3.12}
\left(\begin{matrix} x&y\\ z&w\end{matrix}\right) \longmapsto 
\left(\begin{matrix} 0&y\\ z&0\end{matrix}\right),\qquad x,y,z,w\in \bb M_n,
\end{equation}
is a complete contraction and so  $\widetilde\psi\colon \ \bb M_{2n}\to B(H)$, 
defined by
\begin{equation}\label{eq3.13}
\widetilde\psi \left(\begin{matrix} x&y\\ z&w\end{matrix}\right) = 
(\psi(y) + \psi(z^*)^*)/2,\qquad x,y,z,w\in\bb M_n,
\end{equation}
is a self--adjoint complete contraction. Then
\begin{equation}\label{eq3.14}
\widetilde\psi(\tilde\phi(a)) = (\psi(\phi(a)) + \psi(\phi(a^*))^*)/2,\qquad 
a\in\cl A,
\end{equation}
and \eqref{eq3.10} follows for the self-adjoint elements $x\in\omega$.
\end{proof}

We now show that entropy lifts from a $C^*$-algebra to its unitization.

\begin{mythm}\label{thm3.5}
Let $\cl A$ be a non-unital $C^*$-algebra. If $\alpha$ is an automorphism of 
$\cl A$, then $ht'_{\cl A}(\alpha) = ht'_{\tilde{\cl A}}(\tilde\alpha)$, where 
$\tilde\alpha$ is the unique extension of $\alpha$ to the unitization 
$\tilde{\cl A}$.
\end{mythm}

\begin{proof}
If $\cl A$ is faithfully represented on a Hilbert space $H$, then so too is 
$\tilde{\cl A}$. The inequality $ht'_{\cl A}(\alpha) \le ht'_{\tilde{\cl 
A}}(\tilde\alpha)$ is then trivial. To prove the reverse inequality, we may 
restrict attention to a finite subset $\omega$ of the 
positive open unit ball. The group $\{\alpha^n\}_{n\in\bb Z}$ 
gives an action of $\bb Z$ on $\cl A$, so the crossed product construction of 
the previous section allows us to assume that $\alpha$ is spatially 
implemented. 
Thus we may regard $\tilde\alpha$ as an automorphism of $B(H)$, not just of 
$\tilde{\cl A}$.

Fix $\delta\in (0,1)$ and choose $a\in\cl A^+$, $\|a\| < 1$, such that $x\le 
a$ 
for $x\in\omega$, which is possible since the positive part of the open unit 
ball of any $C^*$-algebra is upward filtering, \cite[Example 3.1.1]{Mur}. 
Define $f_\delta\colon \ [0,1]\to 
[0,1]$ by $f_\delta(t) = t/\delta$ for $t\in [0,\delta]$, and $f_\delta(t) = 
1$ 
elsewhere. Let $f = f_\delta(a) \in\cl A$, and let $p\in B(H)$ be the spectral 
projection of $a$ for the interval $[\delta,1]$. The relations
\begin{equation}\label{eq3.15}
fp=pf = p,\quad ap=pa,\quad \|a-ap\| \le \delta
\end{equation} 
are immediate from the functional calculus. Now enlarge $\omega$ by defining 
$\omega' = \omega\cup \{f,a\}$, and fix an integer  $N$. Then let $n$ denote 
$rcc\left(\bigcup\limits^{N-1}_{i=0} \alpha^i(\omega'),\delta\right)$. By 
definition, 
there exist complete contractions
\begin{equation}\label{eq3.16}
\cl A  \stackrel{\phi}{\longrightarrow} \bb M_n 
\stackrel{\psi}{\longrightarrow} 
B(H)
\end{equation}
such that 
\begin{equation}\label{eq3.17}
\|\psi(\phi(x)) - x\| < \delta,\qquad x\in \bigcup\limits^{N-1}_{i=0} 
\alpha^i(\omega').
\end{equation}
By Lemma \ref{lem3.4}, we may assume that $\phi$ and $\psi$ are self-adjoint. 
Here we may ignore the doubling of dimension which does not affect the 
entropy. 
By injectivity of $\bb M_n$, we may extend $\phi$ to a self-adjoint complete 
contraction on $\tilde{\cl A}$ which we also denote by $\phi$. Then $t = 
\psi(\phi(1)) \in B(H)$ is self-adjoint, so we may define $q$ to be its 
spectral 
projection for the interval $[1-2\delta^{1/2},1]$. The functional calculus 
gives
\begin{equation}\label{eq3.18}
1-t \ge 2\delta^{1/2}(1-q).
\end{equation}
Since $\|1-2f\|\le 1$, it follows that $\|t-2\psi(\phi(f))\| \le 1$, and thus 
\begin{equation}\label{eq3.19}
\|t-2f\| \le \|\psi(\phi(1-2f))\|+2\|\psi(\phi(f)) - f\| \le 1+2\delta.
\end{equation}
We see that $t-2f$ is bounded below by $-(1+2\delta)$, and this inequality,
after rearrangement of the terms, is
\begin{equation}\label{eq3.20}
2f-t-1\le 2\delta.
\end{equation}
We may now multiply on both sides by $p$ to reach
\begin{equation}\label{eq3.21}
p-ptp \le 2\delta p,
\end{equation}
using \eqref{eq3.15}. Then, from \eqref{eq3.18} and \eqref{eq3.21},
\begin{equation}\label{eq3.22}
2\delta p \ge p(1-t)p \ge 2\delta^{1/2} p(1-q)p,
\end{equation}
from which it follows that $\|p(1-q)\| \le \delta^{1/4}$. This yields the 
estimate
\begin{align}
\|a-aq\| &\le \|a-ap\| + \|a(p-pq)\| + \|(ap-a)q\|\nonumber\\
\label{eq3.23}
&\le 2\delta + \delta^{1/4} \le 3\delta^{1/4},
\end{align}
using \eqref{eq3.15}. By Lemma \ref{lem3.1} and \eqref{eq3.23}, if 
$x\in\omega$ 
then
\begin{align}
\|x-qxq\| &\le \|x-qx\| + \|q(x-xq)\|\nonumber\\
&\le 2\|x-xq\|\nonumber\\
&\le 2\|a-aq\|^{1/2}\nonumber\\
\label{eq3.24}
&\le 4\delta^{1/8}.
\end{align}
By repeating this argument with $\alpha^i(a)$, $\alpha^i(f)$ and 
$\tilde\alpha^i(p)$ 
replacing respectively $a,f$ and $p$, we see that \eqref{eq3.24} also holds 
for 
$x\in \bigcup\limits^{N-1}_{i=0} \alpha^i(\omega)$.

Since $\cl A$ is a closed ideal in $\tilde{\cl A}$, we may choose a state 
$\mu$ 
on $\tilde{\cl A}$ which annihilates $\cl A$. Then define $\tilde\phi \colon \ 
\tilde{\cl A}\to \bb M_{n+1}$ by
\begin{equation}\label{eq3.25}
\tilde\phi(x) = \left(\begin{matrix} \phi(x)&0\\ 0&\mu(x)\end{matrix}\right), 
\qquad x\in \tilde{\cl A},
\end{equation}
and define $\widetilde\psi \colon \ \bb M_{n+1}\to B(H)$ to be any completely 
contractive extension  of the map
\begin{equation}\label{eq3.26}
\left(\begin{matrix} m&0\\ 0&\lambda\end{matrix}\right) \longmapsto q\psi(m) q 
+ 
\lambda(1-q),\qquad m\in\bb M_n, \ \lambda\in\bb C.
\end{equation}
If $x\in \bigcup\limits^{N-1}_{i=0} \alpha^i(\omega)$, then \eqref{eq3.24} 
implies that 
\begin{equation}\label{eq3.27}
\|\widetilde\psi(\tilde \phi(x)) - x\| \le 4\delta^{1/8}.
\end{equation}
 From \eqref{eq3.26}
\begin{equation}\label{eq3.28}
\widetilde\psi(\tilde\phi(1)) = qtq + 1-q,
\end{equation}
so by the definition of $q$,
\begin{equation}\label{eq3.29}
\|\widetilde\psi(\tilde\phi(1)) - 1\|\le 2\delta^{1/2}.
\end{equation}
Passage from $\bb M_n$ to $\bb M_{n+1}$ does not affect the entropy, and so we 
have proved that
\begin{equation}\label{eq3.30}
ht'_{\tilde{\cl A}}(\tilde\alpha, \{1\}\cup \omega, 4\delta^{1/8}) \le 
ht'_{\cl A}(\alpha,\omega',\delta),
\end{equation}
and the result follows from Lemma \ref{lem2.1} (ii) with $g(\delta) = 
4\delta^{1/8}$.
\end{proof}

We now turn to the question of whether $ht(\alpha) = ht'(\alpha)$.

\begin{proposition}\label{pro3.6}
Let $\cl A$ be a unital $C^*$-algebra and let $\omega$ be a finite subset of 
$\cl A_1$ containing the identity. For $\delta\in (0,1)$,
\begin{equation}\label{eq3.31}
rcc(\omega,\delta) \ge rcp(\omega,12\delta^{1/8}).
\end{equation}
\end{proposition}

\begin{proof}
Let $n=rcc(\omega,\delta)$ and choose complete contractions
\begin{equation}\label{eq3.32}
\cl A \stackrel{\phi}{\longrightarrow} \bb M_n 
\stackrel{\psi}{\longrightarrow} 
B(H)
\end{equation}
so that
\begin{equation}\label{eq3.33}
\|\psi(\phi(x)) - x\| < \delta,\qquad x\in\omega.
\end{equation}
By polar decomposition followed by diagonalization, there is a non-negative 
diagonal matrix $D$ and two unitaries so that
\begin{equation}\label{eq3.34}
\phi(1) = UDV.
\end{equation}
By replacing $\phi$ by $U^*\phi(\cdot)V^*$ and $\psi$ by $\psi(U\cdot V)$, we 
may assume that $\phi(1) = D$, and it is clear that $\|D\| \in (1-\delta,1]$. 
Let $E$ be the projection onto the space spanned by eigenvectors of $D$ 
corresponding to eigenvalues  in the interval $[1-\delta^{1/2},1]$.

If $X \in (1-E)\bb M_n(1-E)$, $\|X\| \le 1$, then
\begin{equation}\label{eq3.35}
\|D+\delta^{1/2} e^{i\theta}X\| \le 1,\qquad \theta\in \bb R.
\end{equation}
Apply $\psi$ to obtain
\begin{equation}\label{eq3.36}
\|1+\delta^{1/2}e^{i\theta} \psi(X)\| \le \|\psi(D+\delta^{1/2} 
e^{i\theta}X)\| 
+ \delta,\qquad \theta\in\bb R,
\end{equation}
and it follows from Lemma \ref{lem3.2} that
\begin{equation}\label{eq3.37}
\|\psi(X)\| \le 2\delta^{1/2}.
\end{equation}
In particular, \eqref{eq3.37} is valid for $1-E$ and $D(1-E)$. Since $\|DE-E\| 
\le \delta^{1/2}$, we obtain
\begin{align}
\|\psi(D)-\psi(E)\| &\le \|\psi(DE-E)\| + \|\psi(D(1-E))\|\nonumber\\
\label{eq3.38}
&\le 3\delta^{1/2},
\end{align}
from \eqref{eq3.37}. Thus
\begin{equation}\label{eq3.39}
\|1-\psi(E)\| \le \delta+3\delta^{1/2} < 4\delta^{1/2}.
\end{equation}
Using \eqref{eq3.37} with $X = 1-E$, this leads to
\begin{equation}\label{eq3.40}
\|1-\psi(1)\| < 6\delta^{1/2}.
\end{equation}
By Lemma \ref{lem3.3}, there is a completely positive contraction 
$\psi_1\colon 
\ \bb M_n\to B(H)$ so that
\begin{equation}\label{eq3.41}
\|\psi-\psi_1\|_{cb} < 4\delta^{1/4}.
\end{equation}
 From \eqref{eq3.37},
\begin{equation}\label{eq3.42}
\|\psi_1(1-E)\| < 6\delta^{1/4}.
\end{equation}
Let $V^*\pi(\cdot)V$ be the Stinespring representation, \cite{Sti}, of 
$\psi_1$. 
Then \eqref{eq3.42} implies the inequality $\|V^*(1-E)\| < 
6^{1/2}\delta^{1/8}$, so
\begin{equation}\label{eq3.43}
\|\psi_1((1-E)X)\| \le 6^{1/2}\delta^{1/8}\|X\|,\qquad X\in\bb M_n,
\end{equation}
with a similar estimate for $\psi_1(X(1-E))$. Since
\begin{equation}\label{eq3.44}
Y-EYE = Y(1-E) + (1-E)YE,\qquad Y\in\bb M_n,
\end{equation}
we obtain
\begin{equation}\label{eq3.45}
\|\psi_1(Y) - \psi_1(EYE)\| \le 2\cdot6^{1/2} \delta^{1/8}\|Y\|
\end{equation}
from \eqref{eq3.43}. Let $m$ be the rank of $E$ and identify $E\bb M_nE$ with 
$\bb M_m$. Then let $\psi_2\colon \ \bb M_m\to B(H)$ be the restriction of 
$\psi_1$ to $E\bb M_nE$. Also define $\phi_1\colon \ A\to \bb M_m$ by 
$\phi_1(\cdot) = E\phi(\cdot)E$. Then $\|\phi_1(1) - E\| \le \delta^{1/2}$, so 
by Lemma~\ref{lem3.3} there exists a completely positive contraction 
$\phi_2\colon \ A\to \bb M_n$ such that $\|\phi_1-\phi_2\|_{cb} \le 2^{1/2} 
\delta^{1/4}$. Putting together the estimates, we obtain, for $x\in\omega$,
\begin{align}
\|x-\psi_2(\phi_2(x))\| &\le \|x-\psi_2\phi_1(x)\| + 2^{1/2}\delta^{1/4} 
\nonumber\\
&= \|x-\psi(E\phi_1(x)E)\| + 2^{1/2} \delta^{1/4}\nonumber\\
&\le \|x-\psi_1(\phi(x))\| + 2\cdot 6^{1/2} \delta^{1/8} + 2^{1/2} 
\delta^{1/4}\nonumber\\
&\le \|x-\psi(\phi(x))\| + (4+2^{1/2})\delta^{1/4} + 2\cdot6^{1/2} 
\delta^{1/8}\nonumber\\
\label{eq3.46}
&< 12\delta^{1/8},
\end{align}
using \eqref{eq3.45} and \eqref{eq3.33}. This proves \eqref{eq3.31}.
\end{proof}

This completes the technical results needed to prove the following.

\begin{mythm}\label{thm3.7}
Let $\alpha$ be an automorphism of $\cl A$. Then
\begin{equation}\label{eq3.47}
ht(\alpha) = ht'(\alpha).
\end{equation}
\end{mythm}

\begin{proof}
If $\cl A$ is unital then the non-trivial inequality $ht(\alpha) \le 
ht'(\alpha)$ follows from Proposition~\ref{pro3.6}. If $\cl A$ is non-unital 
then we also use Theorem~\ref{thm3.5} to conclude that
\begin{equation}\label{eq3.48}
ht'_{\cl A}(\alpha) = ht'_{\tilde{\cl A}}(\tilde\alpha) = ht_{\tilde{\cl 
A}}(\tilde\alpha) \ge ht_{\cl A}(\alpha),
\end{equation}
proving the result.
\end{proof}

\begin{myrem}\label{rem3.8}
The idea for Proposition \ref{pro3.6} comes from \cite{Sm}, where it appears 
in 
a different context. Henceforth we will refer to entropy  as $ht(\alpha)$, but 
use the definition of $ht'(\alpha)$.$\hfill\square$
\end{myrem}\newpage

\section{Duality}\label{sec4}

\indent

One of the most useful results in the theory of crossed products
is the Imai--Takai duality theorem, \cite{IT}, which generalizes to locally 
compact group actions
the Takai duality theorem, \cite{T}, for abelian groups. This asserts the 
existence of a dual
coaction $\hat\alpha$ on $\AcG$ so that 
\begin{equation}
\coA \approx \AK.
\end{equation}
The new crossed product algebra $\coA$ is formed by combining copies of 
$\AcG$ and $C_0(G)$ acting on $L^2(G,H)\otimes L^2(G)$. In this 
representation, the 
typical generator $\int_G f(s)\tilde\pi(a)\lambda_s\,ds$ of the crossed 
product becomes
\begin{equation}
\int_G f(s)\tilde\pi(a)\lambda_s\otimes l_s\,ds,\ \ f\in C_c(G),\ a\in \cl A,
\end{equation}
where $l_s$ is left translation on $L^2(G)$ (and $r_s$ will denote the 
corresponding
right regular representation). The algebra $C_0(G)$ is represented on the 
second copy of $L^2(G)$
as multiplication operators. Then $\coA$ is the norm closed span of operators 
of the form
\begin{equation}
\int_Gf(s)\tilde\pi(a)\lambda_s\otimes M_Fl_s\,ds,\ \ f \in C_c(G),\ F\in 
C_0(G),\ a\in \cl A.
\end{equation}
In Theorem \ref{thm4.2} we state a modified version of the Imai--Takai duality 
theorem, essentially due to Quigg, \cite{Q}. This is in order to 
obtain information
which is very difficult to extract from the original argument of \cite{IT}. 
 We then give 
some consequences to be used in Section \ref{sec5}.

In \cite{NS}, the notation $f\cdot a$ (where $f \in C_c(G)$ and $a \in \cl A$) 
was used
to denote the function $t\mapsto f(t)a$ in $C_c(G,\cl A)\subseteq L^1(G,\cl 
A)$. The crossed product
is the $C^*$--completion of a representation $\tilde\pi\times \lambda$ of 
$L^1(G,\cl A)$,
where
\begin{equation}
\tilde\pi\times \lambda(f\cdot a)=\int_G f(s)\tilde\pi(a) \lambda_s\,ds.
\end{equation}
Below we will simplify notation and use $f\cdot a$ when we really mean 
its image $\tilde\pi\times \lambda(f\cdot a)
\in \AcG$ under a faithful representation.  It should also be noted that
all maps considered throughout this section are well defined (see for
example ~\cite{NS}).

\begin{mythm}
\label{thm4.2}
  There exists an isomorphism
  \begin{equation}
    \imath:\coA\to\AK
  \end{equation}
which induces a covariant isomorphism between the two $C^*$-dynamical systems 
\begin{equation*}(\coA,\adlrg,G) {\text{ and }}
(\AK,\alpha_g\otimes\adl,G).
\end{equation*}
\end{mythm}
\begin{proof}
  This is a special case of \cite[Theorem 3.1]{Q}, with minor changes.  The 
required isomorphism 
is  given by
   \begin{equation}\label{eq4.15}
    \imath(\tilde\pi(a)\lambda_s\otimes M_Fl_s)=
    (1\otimes M_F)\tilde\pi(a)\lambda_s,
  \end{equation}
for any $a\in \cl A$, $f\in C_c(G)$ and $F\in C_c(G)$.
  \end{proof}

We now begin to construct complete contractions which are covariant with 
respect to 
certain automorphims, and we note that they are variants of maps considered in 
\cite{NS}. Below, $\Delta(\cdot)$ denotes the modular function on a group $G$.

\begin{proposition}
\label{pro4.3}
    Let $g\in G$ satisfy $\Delta(g)=1$ and $f\in C_c(G)$ be such that
    $f(t)=f(gtg^{-1})$ for all $t\in G$. Define the map $S_f:\cl A\to\AcG$
    by $S_f(a)=\int_Gf(s)\tilde\pi(a)\lambda_s\,ds$.

    Then $\|S_f\|_{cb}\le\|f\|_1$ and the following diagram is commutative:
    \begin{equation}
      \xymatrix{\cl A\ar[rr]^{S_f} \ar[d]_{\alpha_g} &&
        \AcG\ar[d]^{Ad_{\lambda_g}}\\ 
        \cl A\ar[rr]_{S_f} && \AcG}
    \end{equation}
  \end{proposition}

  \begin{proof}
    Fix an arbitrary element $a\in \cl A$. Then
    \begin{eqnarray}
      \lambda_gS_f(a)\lambda_{g^{-1}}&=&
      \int_Gf(s)\lambda_g\tilde\pi(a)\lambda_s\lambda_{g^{-1}}\,ds\nonumber\\ 
      &=&\int_Gf(s)\tilde\pi(\alpha_g(a))\lambda_{gsg^{-1}}\,ds\nonumber\\ 
      &=&\int_Gf(g^{-1}sg)\tilde\pi(\alpha_g(a))\lambda_s\,ds\nonumber\\
      &=&S_f(\alpha_g(a)).
    \end{eqnarray}
    The inequality $\|S_f\|_{cb}\le\|f\|_1$ was proved in \cite{NS}.
  \end{proof}

  \begin{proposition}
\label{pro4.4}
    Let $g\in G$ satisfy $\Delta(g)=1$ and suppose that $\eta\in L^2(G)$ is 
such that
    $\|\eta\|_2=1$ and $\eta(gtg^{-1})=\eta(t)$ for almost all $t\in G$.
    Define $T_\eta$ to be the right slice map by $\omega_\eta$ on
    ${B}(L^2(G,H))$, restricted to $\AcG$. 
    Then $\|T_\eta\|_{cb}=1$ and the following diagram is commutative:
    \begin{equation}
      \xymatrix{\AcG\ar[rr]^{T_\eta} \ar[d]_{Ad_{\lambda_g}} &&
      \cl A\ar[d]^{\alpha_g}\\ 
      \AcG\ar[rr]_{T_\eta} && \cl A}
    \end{equation}
  \end{proposition}

  \begin{proof}
    Let $a\in \cl A$ and $f\in C_c(G)$. Then
    \begin{equation}
      T_\eta(f\cdot a)=\int_GF_{f,\eta}(t)\alpha_{t^{-1}}(a)\,dt
    \end{equation}
    where $F_{f,\eta}\in C_c(G)$ is  defined by
    \begin{equation}
      F_{f,\eta}(t)=\int_Gf(s)\eta(s^{-1}t)\overline{\eta(t)}\,ds.
    \end{equation}
    Then
    \begin{equation}
      \alpha_g(T_\eta(f\cdot a))=\int_GF_{f,\eta}(t)\alpha_{gt^{-1}}(a)\,dt
    \end{equation}
    On the other hand, define a function $h\in C_c(G)$ by
    $h(t)=f(g^{-1}tg)$ for all $t\in G$. Then
    \begin{eqnarray}
      \lambda_g(f\cdot a)\lambda_{g^{-1}}&=&
      \int_Gf(s)\tilde\pi(\alpha_g(a))\lambda_{gsg^{-1}}\,ds\nonumber\\
      &=&\int_Gh(s)\tilde\pi(\alpha_g(a))\lambda_s\,ds\nonumber\\
      &=&h\cdot\alpha_g(a).
    \end{eqnarray}
    Furthermore,
    \begin{eqnarray}
 T_\eta(h\cdot\alpha_g(a))&=&\int_GF_{h,\eta}(t)\alpha_{t^{-1}g}(a)\,dt
\nonumber\\
      &=&\int_GF_{h,\eta}(gtg^{-1})\alpha_{gt^{-1}}(a)\,dt
    \end{eqnarray}
    where $F_{h,\eta}$ is defined as before. To conclude the proof it only
    remains to show that, for any $t\in G$,
    $F_{h,\eta}(gtg^{-1})=F_{f,\eta}(t)$. Indeed,
    \begin{eqnarray}
      F_{h,\eta}(gtg^{-1})&=&\int_Gh(s)\eta(s^{-1}gtg^{-1})
      \overline{\eta(gtg^{-1})}\,ds\nonumber\\ 
      &=&\int_Gf(g^{-1}sg)\eta(g^{-1}s^{-1}gt)\overline{\eta(t)}\,ds\nonumber\\
      &=&\int_Gf(s)\eta(s^{-1}t)\overline{\eta(t)}\,ds\nonumber\\
      &=&F_{f,\eta}(t),
    \end{eqnarray}
as required.
  \end{proof}

\begin{proposition}
\label{pro4.5}
  Let $g\in G$.
  Let $F\in C_0(G)$ be such that, for any $t\in G$, $F(gtg^{-1})=F(t)$.
  Define the map $S_F:\AcG\to\coA$ by
  \begin{equation}
    S_F(f\cdot a)=\int_Gf(s)\tilde\pi(a)\lambda_s\otimes M_Fl_s\,ds
  \end{equation}
  for all $f\in C_c(G)$ and $a\in A$.  Then $\|S_F\|_{cb}\le\|F\|_\infty$
  and the following diagram is commutative:
\begin{equation}
  \xymatrix{\AcG\ar[rr]^{S_F} \ar[d]_{\adlg} && \coA\ar[d]^{\adlrg}\\ 
    \AcG\ar[rr]_{S_F} && \coA    }
\end{equation}
\end{proposition}

\begin{proof}
  Fix elements $a\in \cl A$ and $f\in C_c(G)$. Then
  \begin{equation}
    S_F(\lambda_g(f\cdot a)\lambda_{g^{-1}})=
    \int_Gf(s)\tilde\pi(\alpha_g(a))\lambda_{gsg^{-1}} \otimes
    M_Fl_{gsg^{-1}}\,ds
  \end{equation}
  On the other hand,
  \begin{eqnarray}
    \adlrg(S_F(f\cdot a))&=&
    \int_Gf(s)\tilde\pi(\alpha_g(a))\lambda_{gsg^{-1}}\otimes
    l_gr_gM_Fl_sr_{g^{-1}}l_{g^{-1}}\,ds\nonumber\\ 
    &=&\int_Gf(s)\tilde\pi(\alpha_g(a))\lambda_{gsg^{-1}} \otimes
    M_Fl_{gsg^{-1}}\,ds,
  \end{eqnarray}
so these equations prove commutativity of the diagram.
\end{proof}

\begin{proposition}
\label{pro4.6}
  Let $g\in G$ be such that $\Delta(g)=1$.
  Let $\eta\in L^2(G)$, $\|\eta\|_2=1$ be such that, for almost all $t\in G$,
  $\eta(gtg^{-1})=\eta(t)$.  Let $T_\eta:\coA\to\AcG$ be the restriction of
  the right slice map by $\omega_\eta$. Then $\|T_\eta\|_{cb}=1$ and the
  following diagram is commutative:
\begin{equation}
  \xymatrix{\coA\ar[rr]^{T_\eta} \ar[d]_{\adlrg} && \AcG\ar[d]^{\adlg}\\
    \coA\ar[rr]_{T_\eta} && \AcG }
\end{equation}
\end{proposition}

\begin{proof}
  For simplicity put $T=T_\eta$. Let
  \begin{equation}
    y=\int_Gf(s)\tilde\pi(a)\lambda_s\otimes M_Fl_s\,ds
  \end{equation}
  be an element of $\coA$, where $f,F\in C_c(G)$. Then
  \begin{eqnarray}
    \lambda_gT(y)\lambda_{g^{-1}}&=&\int_Gf(s)\omega_\eta(M_Fl_s)
    \tilde\pi(\alpha_g(a))\lambda_{gsg^{-1}}\,ds\nonumber\\
    &=&\int_Gf(g^{-1}sg)\omega_\eta(M_Fl_{g^{-1}sg})\tilde\pi
    (\alpha_g(a))\lambda_s\,ds.
  \end{eqnarray}
  On the other hand let $\tilde F\in C_c(G)$ be defined by $\tilde
  F(t)=F(g^{-1}tg)$ for all $t\in G$. Since
  $l_gr_gM_Fl_sr_{g^{-1}}l_{g^{-1}}=M_{\tilde F}l_{gsg^{-1}}$, we
  also have
  \begin{eqnarray}
    (\lambda_g\otimes l_gr_g)y(\lambda_{g^{-1}}\otimes r_{g^{-1}}l_{g^{-1}}
    )&=&\int_Gf(s)\tilde\pi(\alpha_g(a))\lambda_{gsg^{-1}}\otimes
    M_{\tilde F}l_{gsg^{-1}}\,ds\nonumber\\
    &=&\int_Gf(g^{-1}sg)\tilde\pi(\alpha_g(a))\lambda_s\otimes
    M_{\tilde F}l_s\,ds.
  \end{eqnarray}
  Therefore
  \begin{equation}
    T(\adlrg(y))=\int_Gf(g^{-1}sg)\omega_\eta(M_{\tilde F}l_s)\tilde
    \pi(\alpha_g(a))\lambda_s\, ds.
  \end{equation}
  It remains to observe that, for any $s\in G$, 
  \begin{equation}
    \omega_\eta(M_{\tilde F}l_s)=
  \omega_\eta(M_Fl_{g^{-1}sg}).
  \end{equation}
  Indeed,
  \begin{eqnarray}
    \omega_\eta(M_Fl_{g^{-1}sg})
    &=&\int_G(M_Fl_{g^{-1}sg}\eta)(r)\overline{\eta(r)}\,dr\nonumber\\
    &=&\int_GF(r)\eta(g^{-1}s^{-1}gr)\overline{\eta(r)}\,dr\nonumber\\
    &=&\int_GF(r)\eta(s^{-1}grg^{-1})\overline{\eta(grg^{-1})}\,dr\nonumber\\
    &=&\int_GF(g^{-1}rg)\eta(s^{-1}r)\overline{\eta(r)}\,dr\nonumber\\
    &=&\omega_\eta(M_{\tilde F}l_s),
  \end{eqnarray}
completing the proof.
\end{proof}\newpage

\section{Main results}\label{sec5}

\indent

In this section we will use our preceding results to compare the entropies
of various automorphisms. 
Given a subgroup $H$ of continuous automorphisms of $G$, a function $f$ on
$G$ is $H$-invariant if, for all $\phi\in H$, $f(\phi(x))=f(x)$ everywhere
or almost everywhere, as appropriate. If $H$ is generated by $Ad_g$ for
some fixed $g\in G$, then we may also use the terminology $g$--invariant. In
the case when $H$ is the closure of the inner automorphisms, we will just
say that $f$ is invariant. A set $N$ is called $g$--invariant if $gNg^{-1}=N$.

\begin{lemma}\label{lem5.2}
Let $G$ be a locally compact group and let $g$ be a fixed element of $G$.
\begin{itemize}
\item[(i)]
If the group has a basis  of compact $g$--invariant neighborhoods of $e$, then 
$\Delta(g)=1$;
\item[(ii)]
If the closed subgroup of $Aut(G)$ generated by the inner automorphism $Ad_g$ 
is compact, then there is a basis of compact $g$--invariant neighborhoods of 
$e$.
\end{itemize}
\end{lemma}
\begin{proof}
Any
  compact invariant neighborhood $N$ satisfies $gNg^{-1}=N$, from which the
  equality $\Delta(g)=1$ follows. To prove the second part,
  consider $g \in G$ and let $H$ be the closed subgroup of
  ${\mathrm{Aut}}(G)$ generated by $Ad_g$, which is compact by hypothesis.
  It follows from \cite[Theorem 4.1]{GM} that $G$ has small $H$--invariant
  neighborhoods of the identity, as required.  
\end{proof}

Recall that a group is said to be $[SIN]$ if there is a basis of compact
invariant neighborhoods of $e$. Such groups are unimodular, \cite{mosak1}. If, 
for each $g\in G$, a basis of
$g$--invariant neighborhoods can be found, then we say that $G$ is
{\it {locally [SIN]}}. It is then clear that the following theorem applies
to both classes of groups, while being more general. 
\begin{mythm}\label{thm5.1}
  Let $g\in G$ be a fixed element for which the identity has a basis of
  compact $g$--invariant neighborhoods. Then
    \begin{equation}
      ht_{\cl A}(\alpha_g)\le ht_{\AcG}(\adlg).
    \end{equation}
  \end{mythm}
  \begin{proof}
    We prove that the pair $(\cl A,\alpha_g), (\AcG,\adlg)$ satisfies the
    CCCFP of Section \ref{sec2}.  Given $\varepsilon>0$ and
    $\omega=\{a_1,a_2,\ldots,a_n\}\subseteq \cl A_1$, choose $N$
    to be a compact symmetric $g$-invariant neighborhood of $e$ on which
    \begin{equation}\label{eq5.2}
      \|\alpha_{t^{-1}}(a_i)-a_i\|\le\varepsilon/2,\quad1\le i\le n.
    \end{equation}
    
    Let $\eta$ be the characteristic function of $N$, normalized in
    the $L^2$-norm. Clearly $\eta$ is $g$-invariant. Choose
    $N_1\subseteq N$ a compact symmetric $g$-invariant neighborhood of
    $e$ such that
    \begin{equation*}
      \int_G\eta(t)\eta(s^{-1}t)\,dt>1-\varepsilon/2,\quad s\in N_1.
    \end{equation*}
    Pick $f\in C_c(G)$ to be a positive, $g$-invariant function
    supported on $N_1$ (as in~\cite{mosak1}, this can be achieved by
    choosing $N_2$ to be a compact neighborhood of $e$ such that
    $N_2\cdot N_2\subseteq N_1$ and then letting
    $f=\chi_{N_2}*\tilde\chi_{N_2})$. Upon normalizing $f$ in the
    $L^1$ norm, we may assume that $\|f\|_1=1$. Let
    \begin{equation}
      h(t)=\eta(t)\int_Gf(s)\eta(s^{-1}t)\,ds,\quad t\in G.
    \end{equation}
    Clearly,
    \begin{equation}\label{eq5.4}
      1\ge\int_Gh(t)\,dt=\int_{N_1}f(s)\int_G\eta(t)\eta(s^{-1}t)\,dt\,ds
      \ge(1-\varepsilon/2)\int_Gf(s)\,ds=1-\varepsilon/2.
    \end{equation}
    
    Consider now the complete contractions $S_f$ and $T_\eta$, introduced
    in Propositions~\ref{pro4.3} and \ref{pro4.4} (see also~\cite{NS}). A 
simple calculation
    gives
    \begin{equation}
      T_\eta S_f(a)-a=\int_Nh(t)(\alpha_{t^{-1}}(a)-a)\,dt+\left(
      \int_Gh(t)\,dt-1\right)a. 
    \end{equation}
But then, for $i=1,2,\ldots,n$,
    \begin{equation}
      \|T_\eta S_f(a_i)-a_i\|\le\varepsilon/2+\varepsilon/2=\varepsilon,
    \end{equation}
using (\ref{eq5.2}), (\ref{eq5.4}), and the assumption $\|a_i\|\leq 1$.
By Lemma \ref{lem5.2}, $\Delta(g)=1$, and so
    the theorem now follows by combining Propositions~\ref{pro4.3}
    and~\ref{pro4.4}.
  \end{proof}

We are now able to state and prove the main result of this paper.

Recall from Section \ref{sec1} that the class of locally $[FIA]^-$ groups
in the next result is defined by the requirement that the group of
automorphisms generated by a fixed but arbitrary inner automorphism should
have compact closure in ${\text{Aut}}(G)$.

  \begin{mythm}\label{thm5.3}
    Let $G$ be an amenable group and let $g\in G$ be such that the closure of
the group generated by $Ad_g$ in Aut(G) is compact. Then
    \begin{equation}
      ht_{\AcG}(\adlg)= ht_{\cl A}(\alpha_g).
    \end{equation}
In particular, equality holds for all elements of an amenable locally
$[FIA]^-$ group $G$. 
  \end{mythm}

  \begin{proof}
The inequality $ht_{\AcG}(\adlg)\geq ht_{\cl A}(\alpha_g)$ is an immediate 
consequence of the preceding
three results, and so we consider only the reverse.
    We will prove that the pair $(\AcG,\adlg), (\coA,\adlrg)$ has the CCCFP, 
which will imply 
     that 
    \begin{equation}
      ht_{\AcG}(\adlg)\le ht_{\coA}(\adlrg).
    \end{equation}
    Using the version of the Imai-Takai duality theorem stated in 
Theorem~\ref{thm4.2} and 
basic
    properties of the topological entropy from \cite{Br}, we will obtain
    \begin{eqnarray}
      ht_{\coA}(\adlrg)&=&ht_{\AK}(\alpha_g\otimes\adl)\nonumber\\
      &\le&ht_{\cl A}(\alpha_g)+ht_{\mathcal{K}(L^2(G))}(Ad_{l_gr_g}).
    \end{eqnarray}
    The result will then follow from the vanishing of the topological
    entropy of any automorphism of the algebra of compact operators.
    This is folklore for which we have no reference. However, it is
    simple to prove by applying the definition to finite sets of rank
    one projections, whose span is norm dense in ${\cl K}(H)$.
    
    To prove that the pair $(\AcG,\adlg), (\coA,\adlrg)$ satisfies the CCCFP
    it is enough to consider finite sets of the form
    \begin{equation}
      \omega=\{f_1\cdot a_1,f_2\cdot a_2,\ldots,f_n\cdot a_n\}\subseteq \AcG,
    \end{equation}
    where $a_i\in \cl A$ and $f_i\in C_c(G)$ for $1\le i\le n$. Let $K$ be a
    compact set containing the supports of all functions $f_i$, $1\le i\le
    n$. Fix $\varepsilon>0$, let 
    \begin{equation}
      M=\max_{1\le i\le n}\|f_i\|_1\,\|a_i\|,
    \end{equation}
    and define $\delta$ to be $\varepsilon^2/4M^2$.

    Denote by $H$ the closed subgroup of ${\text{Aut}}(G)$ generated by the 
inner
    automorphism $Ad_g$. Since the set
    \begin{equation}
      \{\phi(t)\,:\,t\in K\cup K^{-1},\phi\in H\}
    \end{equation}
    is compact, we may assume that $K$ is $H$-invariant and $K=K^{-1}$.  By
    amenability of $G$ and using \cite[Proposition 7.3.8]{Ped}, there exists
    $f\in C_c(G)$, $f\ge0$, $\|f\|_1=1$, such that
    \begin{equation}
      \|{}_sf-f\|_1<\delta,\quad s\in K,
    \end{equation}
    where ${}_sf(\cdot)$ is defined by ${}_sf(t)=f(s^{-1}t)$.  As
    in~\cite{mosak2}, define the function $f^\#\in L^1(G)$ by
    \begin{equation}
      f^\#(t)=\int_Hf(\phi(t))\,d\phi,\quad t\in G,
    \end{equation}
    where $d\phi$ is normalized Haar measure on the compact abelian group
    $H$. Then $f^\#\in C_0(G)$, being a vector integral of elements of
    $C_c(G)$ over a compact group.  Since $f\ge0$,
    \begin{align}
      \|f^\#\|_1&=\int_G\int_Hf(\phi(t))\,d\phi\,dt\nonumber\\
      &=\int_H\int_Gf(\phi(t))\,dt\,d\phi
      =\int_H\|f\|_1\,d\phi=1.
    \end{align}
    The last equality holds because any continuous automorphism $\phi$ 
in the closure of the inner automorphisms of
    a unimodular group preserves the Haar measure. Moreover, by invariance
    of $d\phi$,
    \begin{equation}
      f^\#(gtg^{-1})=f^\#(t),\quad t\in G.
    \end{equation}
    Also, for $s\in K$,
    \begin{eqnarray}
      \|{}_sf^\#-f^\#\|_1&=&\int_G\left|
        \int_H\left[f(\phi(s^{-1})\phi(t))-f(\phi(t))\right]\,d\phi\right|
      \,dt\nonumber\\
      &\le&\int_H\int_G\left|{}_{\phi(s)}f(\phi(t))-f(\phi(t))\right|
      \,dt\,d\phi\nonumber\\
      &=&\int_H\|{}_{\phi(s)}f-f\|_1\,d\phi<\delta,
    \end{eqnarray}
    since $\phi(s)\in K$. 
    Let $\eta=\sqrt{f^\#}$. Then $\|\eta\|_2=1$ and by~\cite[p.126]{Pat}, for
    any $s\in K$
    \begin{equation}
      \|{}_s\eta-\eta\|_2\le\|{}_sf^\#-f^\#\|^{1/2}_1<\sqrt{\delta}.
    \end{equation}
    It follows that, for any $s\in K$,
    \begin{equation}
      1-\sqrt{\delta}<\int_G\eta(s^{-1}t)\eta(t)\,dt\le1,
    \end{equation}
    using the Cauchy--Schwarz inequality and the preceding estimate.
    Now choose $L\subseteq G$ to be a large $H$-invariant compact set such that
    \begin{equation}
      \int_L(\eta(t))^2\,dt>1-\sqrt{\delta}.
    \end{equation}
    Then there exists an $H$-invariant function $F\in C_0(G)$  such that
    $F|_L=1$ and $0\le F\le1$ (as before, choose $F_1\in C_0(G)$,
    $F_1|_L=1$, $0\le F_1\le1$ and let $F=F_1^\#$).
    
    To summarize, we have found a positive $g$-invariant function $F\in
    C_0(G)$ and a $g$-invariant unit vector $\eta\in L^2(G)$ such
    that 
    \begin{eqnarray}
      \label{5.22}
      1\ge\int_GF(t)\eta(s^{-1}t)\eta(t)\,dt>1-2\sqrt{\delta},\quad s\in K.
    \end{eqnarray}
    Define the maps $S_F$ and $T_\eta$ as in Propositions~\ref{pro4.5}
    and~\ref{pro4.6} respectively. Then, for any $a\in \cl A$ and
    $f\in C_c(G)$,
    \begin{equation}
      T_\eta S_F(f\cdot a)-f\cdot a=\int_G(h(s)-1)f(s)\tilde\pi(a)
      \lambda_s\, ds,
    \end{equation}
    where $h$ is defined by
    \begin{equation}
      h(s)=\int_GF(t)\eta(s^{-1}t)\eta(t)\,dt.
    \end{equation}
    In particular, for any $1\le i\le n$, 
    \begin{equation}
      \|T_\eta S_F(f_i\cdot a_i)-f_i\cdot a_i\|\le 2\sqrt\delta M=
      \varepsilon,
    \end{equation}
    using the estimate~(\ref{5.22}), which is valid since each $f_i$ is
    supported in $K$.  Since $S_F$ and $T_\eta$ satisfy the conditions of
    Propositions~\ref{pro4.5} and~\ref{pro4.6} respectively, the theorem is
    proved.
  \end{proof}

\end{document}